\documentclass[a4paper,12pt]{article}

\usepackage{amssymb,amsmath}
\usepackage{pstricks}
\date{}
\newtheorem{thm}{Theorem}[section]

\newtheorem{lem}[thm]{Lemma}

\newtheorem{conj}[thm]{Conjecture}

\numberwithin{equation}{section}

\begin{document}
\title{Strong connectivity and directed triangles in oriented graphs. Partial results on a particular case of the Caccetta-H\"{a}ggkvist conjecture}
\author{Nicolas Lichiardopol \\ Lyc\'{e} A. de Craponne, Salon, France\\
e-mail : nicolas.lichiardopol@neuf.fr\\}



\maketitle
\begin{abstract} A particular case of  Caccetta-H\"{a}ggkvist conjecture, says that a digraph of order $n$ with minimum out-degree at least $\frac{1}{3}n$ contains a directed cycle of length at most 3. In a recent paper, Kral, Hladky and Norine (see \cite {Kral}) proved that a digraph of order $n$ with minimum out-degree at least $0.3465n$ contains a directed cycle of length at most 3 (which currently is the best result). A weaker particular case says that a digraph of order $n$ with minimum semi-degree at least $\frac{1}{3}n$ contains a directed triangle. In a recent paper (see \cite{Lichia}), by using the result of \cite{Kral}, the author proved that for $\beta \geq 0.343545$, any digraph $D$ of order $n$ with  minimum semi-degree  at least $\beta n$ contains a directed cycle of length at most $3$ (which currently is the best result).  This means that for a given integer $d\geq 1$, every digraph with minimum semi-degree $d$ and of order $md$ with $m\leq 2.91082$, contains a directed cycle of length at most $3$. In particular, every oriented graph  with minimum semi-degree $d$ and of order $md$ with $m\leq 2.91082$, contains a directed triangle. In this paper, by using again the result of \cite{Kral}, we prove that every oriented graph  with minimum semi-degree $d$, of order $md$ with $2.91082< m\leq 3$ and of strong connectivity at most $0.679d$, contains a directed triangle. This will be implied by  a more general and more precise result, valid not only for $2.91082< m\leq 3$ but also for larger values of $m$. As application, we improve two existing results. The first result (Authors Broersma and Li in \cite{Li}), concerns the number of the directed cycles of length $4$ of a triangle free oriented graph of order $n$ and of minimum semi-degree at least $\frac{n}{3}$. The second result (Authors Kelly, K\"{u}hn and Osthus in \cite{Kuhn}), concerns the diameter of a triangle free oriented graph of order $n$ and of minimum semi-degree at least $\frac{n}{5}$ \end{abstract}
{\it Keywords} : Oriented graph, strong connectivity, girth, triangle
\section{ Introduction and definitions}
\indent The definitions which follow are those of \cite{Bang}.\newline\indent
\indent We consider digraphs without loops and without parallel arcs.
$V(D)$ is the {\it vertex set} of $D$ and
the {\it order} of $D$ is the cardinality of $V(D)$. $\mathcal{A}(D)$ is the set
of the arcs of $D$. We denote by $a(D)$ the number of the arcs of $D$ ({\it size} of $D$).
Two arcs $(x,y)$ and $x',y'$  are {\it independent} if the pairs $\lbrace x,y \rbrace$ and $\lbrace x',y' \rbrace$ are disjoint. \newline \indent We say that a vertex $y$ is  an {\it out-neighbor}  of a vertex $x$
({\it in-neighbour} of $x$) if $(x,y)$ (resp. $(y,x)$) is an arc of $D$. $N_D^+(x)$ is the set of the out-neighbors of $x$ and
$N_D^-(x)$ is the set of the in-neighbors of $x$.
The cardinality of $N_D^+(x)$ is the {\it out-degree} $d_D^+(x)$ of $x$
 and the cardinality of $N_D^-(x)$ is the {\ it in-degree}
$d_D^-(x)$ of $x$. We also put $N_D(x)=N_D^+(x)\cup N_D^-(x)$ and $N'_D(x)=N_D^+(x)\cup N_D^-(x)\cup \lbrace x\rbrace$. When no confusion is possible, we omit the subscript $D$.  We denote by $\delta^+(D)$ the minimum out-degree
of $D$ and by $\delta^-(D)$ the minimum in-degree of $D$. The {\it minimum
semi-degree} of $D$ is $\delta^0(D)=\min
\lbrace\delta^+(D),\delta^-(D\rbrace$.\newline \indent For a vertex $x$ of $D$ and for a subset $S$ of $V(D)$, $N_S^+(x)$  is the set of
the out-neighbors of $x$ which are in $S$, and $d_S^+(x)$ is the cardinality of $N_S^+(x)$. Similarly, $N_S^-(x)$) is the set of
the in-neighbors of $x$ which are in $S$, and $d_S^-(x)$ is the cardinality of $N_S^-(x)$.

\vspace{0.2cm} A {\it directed path} of length $p$ of
$D$ is a list $x_0, \ldots, x_p $ of distinct vertices  such that
$(x_{i-1}, x_i)\in \mathcal{A}(D)$ for $1\leq i\leq p$.
A {\it
directed cycle} of length $p\geq 2$ is a list ($x_0 \ldots, x_{p-1}, x_0) $
of vertices with $x_0 \ldots, x_{p-1}$ distinct, $(x_{i-1},
x_i)\in \mathcal{A}(D)$ for $1\leq i\leq p-1$ and $(x_{p-1}, x_0)\in
\mathcal{A}(D)$. From now on, we omit the adjective " directed". A $p$-cycle of $D$ is a directed cycle of length $p$.
\newline A {\it digon} is a  $2$-cycle, and a
triangle is a  $3$-cycle of $D$ of length $3$. The {\it girth} $g(D)$ of $D$ is the minimum length of the
cycles of $D$. \vspace{0.2cm} The digraph $D$ is said to be {\it strongly connected} (for briefly strong) if for every distinct vertices $x$  and $y$ of $D$, there exists a path from $x$ to $y$. It is known that in a non-strong digraph $D$, there exists a partition $(A,B)$ of $V(D)$ with $A\neq \emptyset$ and $B\neq \emptyset$ such that there are no arcs from a vertex of $B$ to a vertex of $A$. (one say that $A$ dominates $B$). We say that a subset $S$ of $V(D)$ disconnects $D$, if the digraph $D-S$ is non-strong. The {\it strong connectivity} $k(D)$ of $D$ is the smallest of the positive integers $m$ such that there exists a subset of $V(D)$ of cardinality $m$ disconnecting $D$. $D$ is said to be {\it $p$-strong connected} if $k(D)\geq p$. It is well known that in a $p$-strong connected digraph, if $S$ is a subset of $V(D)$ such that $\lvert S \rvert \geq p$ and $\lvert V(D)\setminus S \rvert \geq p$, then there exist $p$ independent arcs with starting vertices in $S$ and with ending vertices in  $V(D)\setminus S$. \newline \indent In a strong digraph $D$, for  vertices $x$ and $y$ of $D$, the {\it distance} $d(x, y)$ from $x$ to $y$ is the length of a shortest path from $x$ to $y$. The {\it diameter} $\mathrm{diam}(D)$ is the maximum of the distances  $d(x, y)$. The {\it eccentricity} $\mathrm{ecc}(x)$ of a vertex $x$ is the maximum of the distances $d(x,y)$, $y\in V(D)$. It is clear that $\mathrm{ecc}(x)\leq \mathrm{diam}(D)$ for every vertex $x$ of $D$.

\vspace{0.2cm}
An {\it oriented graph}, is a digraph $D$ such that for any two
distinct vertices $x$ and $y$ of $D$, at most one of the ordered pairs $(x,y)$
and $(y,x)$ is an arc of $D$. The author proved in \cite{Lichiard} that the strong connectivity $k$ of an oriented graph $D$ of order $n$, satisfy $k\geq \dfrac{2(\delta^+(D)+\delta^-(D)+1)-n}{3}$, and this shows that an oriented graph of order $n$ and of minimum semi-degree at least $\frac{n}{4}$, is strongly connected.

\vspace{0.3cm} Caccetta and H\"{a}ggkvist (see
\cite{Caccetta})  conjectured in 1978 that the girth of any
digraph of order $n$ and of minimum out-degree at least $d$ is at
most $\lceil n/d\rceil$.\newline The conjecture is still open when
$d\geq n/3$, in other words it is not known if any digraph of
order $n$  and  minimum out-degree at least $n/3$ contains a
cycle of length at most 3.
\newline In fact it is also unknown if any digraph of order $n$
with both minimum out-degree and minimum in-degree at least $n/3$
contains a cycle of length at most $3$ and then a special case of the Caccetta-H\"{a}ggkvist conjecture is :
\begin{conj} Every digraph of order $n$ and of minimum semi-degree at least $\frac{n}{3}$, contains a cycle of length at most $3$. \end{conj}

\vspace {0.3cm} Two questions
were naturally raised :\\[0.2cm]
 {\bf Question Q$_1$} What is the minimum constant
$c$ such that any digraph of order $n$ with minimum out-degree at
least $cn$ contains a cycle of length at most
$3$. \\[0.1cm] {\bf Question Q$_2$} What is the minimum constant
$c'$ such that any digraph of order $n$ with both minimum out-degree
and minimum in-degree at least $c'n$ contains a cycle of length at
most $3$.

\vspace{0.2cm} It is known that $c\geq c'\geq 1/3$ and the
conjecture is that $c=c'=1/3$.
In a very recent paper (See \cite{Kral}), Hladk\'{y}, Kr\'{a}l'
and Norine proved that $c\leq 0.3465$, which currently is the best result. \newline  By
using this result, the author proved in \cite{Lichia}
that $c'\leq 0.343545$, which currently is the best result. In other terms, this means :
\begin{thm}  For $d\geq 1$, any digraph with minimum semi-degree $d$ and of order at most $2.91082 d$
contains a cycle of length at most $3$. \end{thm}
In our paper, we will see that in an oriented graph $D$ of minimum semi-degree $d$ and of order $md$ with $2.91082 <m< \frac{2}{c}$, an adequate  upper bound on the connectivity of $D$ forces the existence of a triangle. More precisely, we prove :
\begin{thm} Let $D$ be an oriented graph of minimum semi-degree $d$, of order $n=md$ with $2.91082<m<\frac{2}{c}$.
If the connectivity $k$ of $D$ verifies $k\leq \max \left\lbrace \dfrac{5-m-4c+c^2}{(1-c)(2-c)}d,\,\dfrac{2-c m}{2-c}d\right\rbrace $, then $D$ contains at least a triangle.\end{thm}
Since $c\leq 0.3465$, an easy consequence will be :

\begin{thm} Let $D$ be an oriented graph of minimum semi-degree $d$, of order $n=md$ with $2.91082<m\leq 3$.
If the connectivity $k$ of $D$ verifies $k\leq 0.679 d$, then $D$ contains at least a triangle.\end{thm}

Broersma and Li proved in \cite{Li} that in a triangle-free oriented graph of order $n$ and of minimum semi-degree at least $\frac {n}{3}$, every vertex is in more than $1+\frac{n}{15}(11-4\sqrt{6})$ $4$-cycles. We improve this result by proving :
\begin{thm} Let $D$ be a triangle-free oriented graph of minimum semi-degree $d$, of order $n=md$ with $m\leq 3$. Then every vertex $x$ of $D$ is contained in more than $\dfrac{2(5-m-4c+c^2)d}{(1-c)(2-c)}+(2-m)d+1$ cycles such that two of these cycles have only the vertex $x$ in common.\end{thm} If we allow distinct $4$-cycles with others vertices than $x$ in common, we give an even more spectacular improvement, by proving :
\begin{thm} Let $D$ be a triangle-free oriented graph of minimum semi-degree $d$, of order $n=md$ with $m\leq 3$.\newline Then every vertex $x$ of $D$ is contained in more than $\dfrac{11-15c+7c^2-c^3-(c^2-3c+3)m}{(1-c)^2(2-c)}d$  $4$-cycles. \end{thm}

Kelly, K\"{u}hn and Osthus proved in \cite{Kuhn} that if $D$ is an oriented graph of order $n$ and of minimum semi-degree greater than $\frac{n}{5}$, then either the diameter of $D$ is at most $50$ or $D$ contains a triangle. We will considerably improve this result by proving :
\begin{thm} If $D$ is a triangle-free oriented graph of minimum semi-degree $d$ and of order $n=md$ with $m\leq 5$, then the diameter of $D$ is at most $9$.\end{thm}
A result of Chudnovsky, Seymour and
Sullivan (see\cite{Chud}) asserts that one can delete $k$ edges
from a triangle-free digraph $D$ with at most $k$ non-edges to
make it acyclic. Hamburger, Haxell, and Kostochka used this to
prove in \cite{Hamb} that in a triangle-free digraph $D$ with at
most $k$ non-edges,  $\delta^+(D)<\sqrt {2k}$ (and
$\delta^-(D)<\sqrt {2k}$ also) .
\newline Chen, Karson, and Shen improved in \cite{Chen} the initial result of
\cite{Chud} by asserting that one can delete $0.8616k$ edges from
a triangle-free digraph $D$ with at most $k$ non-edges to make it
acyclic. From this result, by using the reasoning of Hamburger,
Haxell and Kostochka  in \cite{Hamb}, it is easy to prove that in
a triangle-free digraph $D$ with at most $k$ non-edges,
$\delta^+(D)<\sqrt {1.7232k}$ and $\delta^-(D)<\sqrt {1.7232k}$.
As the maximum size of an oriented graph of order $n$ is
$\frac{n(n-1)}{2}$, an immediate consequence is :
\begin{lem} If $D$ is a triangle-free oriented graph of order $n$, then
$a(D)< \dfrac{n^2}{2}-\dfrac{(\delta^+(D))^2}{1.7232}$ and $a(D)<
\dfrac{n^2}{2}-\dfrac{(\delta^-(D))^2}{1.7232}$. \end{lem}

\section{ Proofs of  Theorems 1.3 and 1.4}
By hypothesis, $D$ is an oriented graph of minimum semi-degree $d$, of order $n=md$ with $2.91082<m<\frac{2}{c}$ and of strong connectivity $k$ We put $k'=\frac{k}{d}$. Let $K$ be a set of $k$ vertices disconnecting $D$. Then there exists a partition of $V(D)\setminus K$ into two subsets $A$ and $B$, such that there are no arcs from a vertex of $B$ to a vertex of $A$. Without loss of generality, we may suppose that $\lvert B\rvert \leq \lvert A\rvert$. We put $a =\dfrac{\lvert A\rvert}{d}$ and $b =\dfrac{\lvert B\rvert}{d}$. Since $b\leq a$, it holds $b\leq \dfrac{m-k'}{2}$.   First we claim that :
\begin{lem} If $D$ is triangle-free, then for every arc $(y,x)$ of $D$ with $y\in A$ and $x\in B$, it holds $d_B^+(x)+d_A^-(y)\geq 2d-k'd$. \end{lem}
{\bf Proof.} Since $x$ has no out-neighbors in $A$, $x$ has $d^+(x)-d_B^+(x)$ out-neighbors in $K$, which means $\lvert N^+_K(x)\rvert =d^+(x)-d_B^+(x)$.
Since $y$ has no in-neighbors in $B$, $y$ has $d^-(y)-d_A^-(y)$ in-neighbors in $K$, which means $\lvert N^-_K(y)\rvert =d^-(y)-d_A^-(y)$. Since $N^+_K(x)$
and $N^-_K(y)$ are vertex-disjoint (for otherwise, we would have a triangle), we have $d^+(x)-d_B^+(x)+d^-(y)-d_A^-(y)\leq k'd$, hence $d_B^+(x)+d_A^-(y)\geq d^+(x)+d^-(y)-k'd$ and since $d^+(x)\geq d$ and $d^-(y)\geq d$, the result follows \hspace {0.3cm} $\blacksquare$

\vspace{0.2cm} Now, we claim :
\begin{lem} Suppose that $2.91082<m<5-4c+c^2$.
If the connectivity $k$ of $D$ verifies $k\leq \dfrac{5-m-4c+c^2}{(1-c)(2-c)}d$, then $D$ contains at least a triangle. \end{lem}{\bf Proof.} We put $k'=\frac{k}{d}$. Suppose, for the sake of a contradiction, that $D$ does not contain triangles. Let $sd$ be the minimum out-degree of $D\lbrack B \rbrack$, and let $x$ be a vertex of $B$ with $d_B^+(x)=sd$. It is easy to verify that $\dfrac{5-m-4c+c^2}{(1-c)(2-c)}<1$ and since all the out-neighbors of $x$ are in $B\cup K$, it follows that $N_B^+(x)\neq \emptyset$, and so $s>0$. There exists a vertex $x'$ of $N_B^+(x)$, such that $d_{N_B^+(x)}^+(x')<c sd$. It follows that $x'$ has more than $(s-c s)d=(1-c)sd$ out-neighbors in $B$ but not in $N_B^+(x)$, and these out-neighbors cannot be in-neighbors of $x$ (for otherwise, we would have a triangle). We get then  $d_{B\cup K}^-(x)<\lbrack b+k'-1-(1-c)s\rbrack d$. Suppose that $b+k'-1\geq 1$. Then $k'\geq 2-b$, and since $b\leq \dfrac{m-k'}{2}$, we get $k'\geq 2-\dfrac{m-k'}{2}$, hence $k'\geq 4-m$.
Then, since $k'\leq \dfrac{5-m-4c+c^2}{(1-c)(2-c)}$, we get $4-m \leq \dfrac{5-m-4c+c^2}{(1-c)(2-c)}$, hence $(4-m)(c^2-3c+2)\leq 5-m-4c+c^2$. This yields $m(c^2-3c+1)\geq 3c^2-8c+3$, hence $m(c^2-3c+1)\geq 3(c^2-3c+1)+c$. Since $c^2-3c+1>0$, we get $m\geq 3+ \dfrac{c}{c^2-3c+1}$. It is easy to verify that for $\frac{1}{3}\leq c\leq 0.3465$, it holds $\dfrac{c}{c^2-3c+1}>1$. We get then $m>4$, and it is easy to verify that this is contradictory with $m<5-4c+c^2$. Consequently, we have $b+k'-1<1$. We deduce $d_{B\cup K}^-(x)<d$, which means that $N_A^-(x)\neq \emptyset$ (in fact, by the above reasoning, this is true for every vertex of $B$). More precisely, we have \begin{center}$d_A^-(x)>\lbrack 2-k'-b+(1-c)s \rbrack d $ \hspace{2cm} (1) \end{center}
There exists a vertex $y$ of $N_A^-(x)$ with fewer than $c d_A^-(x)$ in-neighbors in $N_A^-(x)$ (for otherwise $D\lbrack N_A^-(x)\rbrack$ would contain a triangle). It follows $d_A^-(y)<c d_A^-(x)+ad-d_A^-(x)$, hence
$d_A^-(y)<ad-(1-c)d_A^-(x)$. From Lemma 2.1, we get $d_A^-(y)\geq (2-k')d-d_B^+(x)$, that is $d_A^-(y)\geq (2-k'-s)d$. We deduce $(2-k'-s)d<ad-(1-c)d_A^-(x)$, hence \begin{center} $sd>(2-k'-a)d+(1-c)d_A^-(x)$ \hspace{2cm} (2) \end{center}
From (1) and (2), we deduce $sd>(2-k'-a)d+(1-c)\lbrack 2-k'-b+(1-c)s \rbrack d$, hence $s>2-k'-a+2-2c -k'+c k'-b+bc+(1-c)^2s$.
It follows  $(2c-c^2)s>4-2k'-a-b-2c+c k'+bc$, and since $a+b=m-k'$, we get $(2c-c^2)s>4-m-k'-2c+c k'+bc$. Since $s<bc $ (for otherwise $D\lbrack B\rbrack$ would contain a triangle), we get $(2c-c^2)bc>4-m-k'-2c+c k'+bc$, hence
$(1-c)^2bc < m+2c-4+(1-c)k'$.
Since all the out-neighbors of $x$ are in $B\cup K$, we have $1-s\leq k'$, hence $s\geq 1-k'$, and since $s< bc$, we get $bc>1-k'$. It follows
$(1-k')(1-c)^2< m+2c-4+(1-c)k'$, hence $k'(1-c)(2-c)>1-2c+c^2-m-2c+4$. This implies $k'>\dfrac{5-m-4c+c^2}{(1-c)(2-c)}$, which is contradictory with the hypothesis on $k$. Consequently $D$ contains at least a triangle, and so, the result is proved. \hspace {0.3cm} $\blacksquare$

\vspace{0.2cm} We claim also :
\begin{lem} Suppose that $2.91082<m<\frac{2}{c}$.
If the connectivity $k$ of $D$ verifies $k\leq \dfrac{2-c m}{2-c}d$, then $D$ contains at least a triangle.\end{lem}{\bf Proof.} Suppose, for the sake of a contradiction, that $D$ does not contain triangles. Let $sd$ be the minimum out-degree of $D\lbrack B \rbrack$, and let $x$ be a vertex of $B$ with $d_B^+(x)=sd$. We have then $k'\geq 1-s$, hence $s\geq 1-k'$. Since $s<bc$ (for otherwise we would have a triangle), we get $bc>1-k'$. Since $b\leq \dfrac{m-k'}{2}$, it follows $\dfrac{(m-k')c}{2}>1-k'$, hence
$mc-k'c >2-2k'$. It follows $k'> \dfrac{2-c m}{2-c}$, which is contradictory with the hypothesis on $k=k'd$. So, the result is proved. \hspace{0.3cm} $\blacksquare$

\vspace{0.2cm} It is easy to prove that $5-4c+c^2<\frac{2}{c}$. By using these two lemmas, we get Theorem 1.3.

\vspace{0.2cm}
It is easy to see that we have $ \dfrac{5-m-4c+c^2}{(1-c)(2-c)}\geq\dfrac{2-c m}{2-c}$ if and only if $m\leq \dfrac{3-2c+c^2}{1-c+c^2}$. Then Theorem 1.3 means that when $2.91082 < m\leq \dfrac{3-2c+c^2}{1-c+c^2}$, a strong connectivity not greater than $\dfrac{5-m-4c+c^2}{(1-c)(2-c)}d$ forces a triangle in $D$, and when $\dfrac{3-2c+c^2}{1-c+c^2}<m<\frac{2}{c}$, a strong connectivity not greater than $\dfrac{2-c m}{2-c}d$ forces a triangle in $D$.\newline \indent It is easy to see that for $2.91082 < m\leq 3$, we have $m< \dfrac{3-2c+c^2}{1-c+c^2}$. Since $c\leq 0.3465$, it is easy to see that we have $0.679d<\dfrac{5-m-4c+c^2}{(1-c)(2-c)}d$. Then by Lemma 2.2, a strong connectivity no greater than $0.679d$ forces a triangle, and so Theorem 1.4 is proved.
 Since a digraph which is not oriented contains a digon, it is easy to see that  proving  Conjecture 1.1, amounts to proving that every oriented graph, of minimum semi-degree at least $d$, of order $md$ with $2.91082 < m\leq 3$ and of connectivity $k>0.679d$, contains at least a triangle.
\section{ Proofs of  Theorems 1.5, 1.6 and 1.7}
{\bf a) Proof of Theorem 1.5} \\[0.2cm] By hypothesis $D$ is a triangle-free oriented graph of minimum semi-degree $d$, of order $n=md$ with $m\leq 3$, and $x$ is a vertex of $D$. Let $k$ be the strong connectivity of $D$ (and $k'=\frac{k}{d}$). We have $k>0$ (for otherwise, by Theorem 1.3 we would have triangles). Clearly, we have $d^+(x)+d^-(x)<md$, and since $k\leq d^-(x)$, it follows $d^+(x)+k<md$, hence $md - d^+(x)>k$. As we have also $d^+(x)\geq k$, there exist $k$ independent arcs $(y_1,z_1), \ldots, (y_k, z_k)$ with $y_i\in N^+(x)$, $z_i \notin N^+(x)$ and $z_i\neq x$ for $1\leq i\leq k$. Since $D$ is triangle-free, we have also $z_i \notin N^-(x)$ for $1\leq i\leq k$. It follows that the set $S_1=\lbrace z_1,\ldots, z_k\rbrace$ is contained in $V(D)\setminus N'(x)$. Similarly, there exist $k$ independent arcs $(v_1,u_1), \ldots, (v_k, u_k)$ with $u_i\in N^-(x)$, $v_i \notin N^-(x)$ and $v_i\neq x$ for $1\leq i\leq k$. Since $D$ is triangle-free, we have also $v_i \notin N^+(x)$ for $1\leq i\leq k$. It follows that the set $S_2=\lbrace v_1,\ldots, v_k\rbrace$ is contained in $V(D)\setminus N'(x)$. We have $\lvert S_1\cap S_2\rvert =\lvert S_1 \rvert + \lvert S_2 \rvert -\lvert S_1\cup S_2\rvert$. Since $\lvert S_1 \rvert = \lvert S_2 \rvert =k'd$ and $\lvert S_1\cup S_2\rvert$ is contained in $V(D)\setminus N'(x)$, it follows $\lvert S_1\cap S_2\rvert\geq 2k'd-(md-d^+(x)-d^-(x)-1)$, hence $\lvert S_1\cap S_2\rvert\geq 2k'd-md+d^+(x)+d^-(x)+1$. Since $d^+(x)\geq d$ and $d^-(x)\geq d$, it follows $\lvert S_1\cap S_2\rvert\geq (2k'+2-m)d+1$. This implies the existence of at least $(2k'+2-m)d+1$ $4$-cycles containing $x$ and such that any two of these cycles have only $x$ in common. Now since $D$ is triangle-free, we deduce from Theorem 1.3 that $k'>\dfrac{5-m-4c+c^2}{(1-c)(2-c)}$, and then Theorem 1.5 is proved.

\vspace{0.5cm} Since $c\leq 0.3465$ and $m\leq 3$, it is easy to see that the number $n_D(x,4)$  of $4$-cycles of $D$ containing $x$, and such that any two of these cycles have only $x$ in common, is at least $ \dfrac{2\times(5-3-4\times 0.3465+0.3465^2)d}{0.6535\times 1.6535}-d+1$, hence $n_D(x,4)>0.358d+1$, and since $d\geq \frac{n}{3}$ ($n$ being the order of $D$), we get $n_D(x,4)>0.119n+1$. Since $1+\frac{n}{15}(11-4\sqrt{6}) \approx 1+0.08014n$ (exceeding value), it is clear that our result improve that of Broersma and Li.\\[0.4cm]
{\bf b) Proof of Theorem 1.6}

\vspace{0.1cm} Let $k=k'd$ be the strong connectivity of $D$. By Theorem 1.4, we have $k>0.679d$. Clearly the eccentricity ecc$(x)$ of $x$ is at least $3$ (for otherwise,we would have a triangle). The author proved in \cite{Lichiard} that the diameter of an oriented graph of order $n$ and of minimum semi-degree at least $\frac{n}{3}$ is at most $4$. By this result, we have $\mathrm{ecc}(x)\leq 4$, and consequently $3\leq \mathrm{ecc}(x)\leq 4$. For $1\leq i\leq \mathrm{ecc}(x)$ let $R_i$ be the set of the vertices $z$ of $D$ such that $d(x,z)=i$. Since $D$ is triangle-free, all the in-neighbors of $x$ are in $R_3 \cup \cdots\cup R_{\mathrm{ecc}(x)}$.\newline \indent
We claim that $d^-_{R_3}(x)>d-\dfrac{m-2-k'}{1-c}d$ (Assertion (Ass)).
\newline \indent We observe first that $m-2-k'>0$. Indeed, for an arbitrary vertex $u$ of $D$, there exists $k'd$ independent arcs with starting vertices in $N^+(u)$ and ending vertices in $V(D)\setminus N^+(u)$. Since $D$ is triangle-free  these ending vertices are not in $N^-(u)$. It follows $2d+k'd<md$, hence  $m-2-k'>0$.\newline\indent
Suppose first that $\mathrm{ecc}(x)=3$. Then all the in-neighbors of $x$ are in $R_3$. This implies $d^-_{R_3}(x)\geq d$, and since $d > d-\dfrac{m-2-k'}{1-c}d$, the assertion (Ass) is proved. \newline Suppose now that $\mathrm{ecc}(x)=4$. Since $R_2$ disconnects $D$, we have $r_2\geq k'd$. Suppose first that $r_3\geq d$. We have $r_4=md-r_1-r_2-r_3-1$, hence $r_4<md-d-k'd-d$, that is $r_4<(m-2-k')d$. It follows $d^-_{R_3}(x)> d-(m-2-k')d$, and since
$d-(m-2-k')d >d-\dfrac{m-2-k'}{1-c}d$, the Assertion (Ass) is proved. Suppose now that $r_3<d$. Clearly, all the in-neighbors of a  vertex of $R_4$ are in $R_3\cup R_4$. It follows that every vertex of $R_4$ has at least $d-r_3$ in-neighbors in $R_4$. Since $D\lbrack R_3\rbrack$ is triangle-free, it holds $d-r_3<cr_4$, hence $r_4>\dfrac{d-r_3}{c}$, hence $r_4>\dfrac{d-(md-r_1-r_2-r_4)}{c}$, which gives $r_4>\dfrac{(1-m)d+r_1+r_2+r_4}{c}$. Since $r_1\geq d$ and $r_2\geq k'd$, we get $r_4>\dfrac{(2-m+k')d+r_4}{c}$, hence $(1-c)r_4<(m-2-k')d$, and then $r_4<\dfrac{m-2-k'}{1-c}d$. It follows $d^-_{R_3}(x)> d-\dfrac {m-2-k'}{1-c}d$, which is the assertion (Ass). It is easy to see that an  in-neighbor $z$ of $x$ which is in $R_3$ has an in-neighbor $z_2$ in $R_2$ and that $z_2$ has an in-neighbor $z_1$ in $R_1$. Then $C_z=(x, z_1, z_2, z, x)$ is a $4$-cycle of $D$, containing $x$ . It is clear that the cycles $C_z$, $z\in N^-_{R_3}(x)$ are distinct. Consequently the vertex $x$ is contained in more than $d-\dfrac {m-2-k'}{1-c}d$ $4$-cycles. Since $k>\dfrac{5-m-4c+c^2}{(1-c)(2-c)}d$ (By Theorem 1.3), the result follows.

\vspace{0.2cm} Since $c\leq 0.3465$, $m\leq 3$ and $k'>0.679$, it holds $d_{R_3}^-(x)>d-\dfrac{3-2-0.679}{1-0.3465}d$, hence $d_{R_3}^-(x)>0.5087d$, hence
$d_{R_3}^-(x)>0.169 n$. So $D$ possess  more than $0.169 n$ $4$-cycles containing $x$, which is much better that the result of Broersma and Li.
\\[0.4cm]
{\bf c) Proof of Theorem 1.7}

\vspace{0.1cm} By hypothesis $D$ is a triangle-free oriented graph of minimum semi-degree $d$, of order $n=md$ with $m\leq 5$. Suppose, for the sake of a contradiction, that the diameter of $D$ is at least $10$. Then let $x$ and $y$ be two vertices of $D$ such that $d(x,y)\geq 10$. For $1\leq i\leq 6$, let $R_i$ be the set of the vertices $z$ of $D$ such that $d(x,z)=i$, and for $1\leq i\leq 3$, let $R_{-i}$ be the set of the vertices $z$ of $D$ such that $d(z,y)=i$. For $1\leq i \leq 6$, $r_i$ is the cardinality of $R_i$ and for $1\leq i\leq 3$, $r_{-i}$ is the cardinality of $R_{-i}$. The sets $R_i$, $1\leq i\leq 6$ are mutually vertex-disjoint, the sets $R_{-i}$, $1\leq i\leq 3$ are also mutually vertex-disjoint, and a set $R_i$, $1\leq i\leq 6$ is a vertex-disjoint with a set $R_{-j}$, $1\leq j\leq 3$ ( for otherwise the diameter of $D$ would be at most $9$). For $2\leq i\leq 6$ we put $R^{'}_i=R_1\cup \cdots R_i$, for $2\leq i\leq 3$ we put $R^{'}_{-i}=R_{-1}\cup \cdots R_{-i}$, and $r^{'}_i$, $r^{'}_{-i}$ are the respective cardinalities. \newline  We claim that $r'_3\geq 2.239d$. Indeed, since $D\lbrack R_1\rbrack$ is triangle-free, there exists a vertex $u$ of $R_1$ with fewer than $0.3465d$ out-neighbors in $R_1$, and then we have $r_2>0.6535d$, hence $r_1+r_2>1.6535 d$. Now, if $r_3\geq d$, it follows $r'_3\geq 2.6535d$, and the assertion is proved. Suppose now that $r_3<d$. It is easy to see that a vertex of $R_2$ has all its out-neighbors in $R'_3$. It follows that a vertex of $R_2$ has at least $d-r_3$ out-neighbors in $R'_2$. Since every vertex of $R_1$ has all its out-neighbors in $R'_2$, it follows $a(D\lbrack R'_2 \rbrack) \geq r_1d+r_2(d-r_3)$, hence : \begin{center}$a(D\lbrack R'_2 \rbrack) \geq r_1d+r_2d-r_2r_3 $ \hspace{2cm} (3) \end{center}
On the other hand by Theorem 1.7, we have \begin{center}$a(D\lbrack R'_2 \rbrack) \leq \dfrac{(r'_2)^2}{2}-\dfrac{(d-r_3)^2}{1.7232}$ \hspace{2cm} (4) \end{center}

\vspace{0.2cm} From (3) and (4), we deduce $r_1d+r_2d-r_2r_3 \leq \dfrac{r_1^2+r_2^2+2r_1r_2}{2}-\dfrac{d^2-2dr_3+r_3^2}{1.7232}$, hence
$3.4464r_1d+3.4464r_2d-3.4464r_2r_3 \leq 1.7232r_1^2+3.4464r_1r_2+1.7232r_2^2-2d^2+4r_3d-2r_3^2$. An easy calculation yields : $1.7232(r_2+r_3+r_1-d)^2 \geq 3.7232r_3^2-(7.4464d-3.4464r_1)r_3+3.7232d^2$. Since $r_1\geq d$, we get $1.7232(r_2+r_3+r_1-d)^2 \geq 3.7232r_3^2-4r_3d+3.7232d^2$, that is  $1.7232(r_2+r_3+r_1-d)^2 \geq f(r_3)$, $f$ being the function defined by $f(t)=3.7232t^2-4dt+3.7232d^2$. By a classical result on the functions of second degree, we have $f(r_3)\geq f\left (\frac{2d}{3.7232}\right )$, hence $f(r_3)> 2.648d^2$. We deduce then $1.7232(r_2+r_3+r_1-d)^2>2.648d^2$, hence $r_2+r_3+r_1-d>1.239d$ which yields $r'_3>2.239d$, and the assertion is still proved. Similarly, we have $r'_{-3}>2.239$. Since $D$ is triangle-free, by Theorem 1.3, the strong connectivity $k$ of $D$ verifies $k>\dfrac{2-5c}{2-c}d$, and since $c\leq 0.3465$, we get $k>0.161d$. It is clear that each of the sets $R_4$, $R_5$ and $R_6$ disconnects $D$, and then $r_i>0.161d$ for $4\leq i\leq 6$. Suppose that $r_4<0.205d$. Then $D\lbrack R'_3\rbrack$, which is triangle-free, is of minimum out degree at least $0.795d$. It follows $0.795<0.3465 r'_3$, hence $r'_3 > 2.2943d$. We have then $v(D)>2.2943d+2.239d+3\times 0.161d$, that is $v(D)>5.0163d$, which is not possible. It follows $r_4\geq 0.205d$. We deduce then $v(D)>2.239d+2.239d+ 0.205d+2\times 0.161d$, that is $v(D)>5.005d$, which is still impossible. Consequently, the diameter of $D$ is at most $9$, and the result is proved. \hspace{0.3cm} $\blacksquare$
\section{ An open problem} Theorem 1.3 gives rise to the following question :\\[0.3cm]
{\bf Open Problem }. For $r$ with $2<r< \frac{2}{c}$, what is the maximum number $\psi (r)\in \rbrack 0, 1\rbrack$ such that every oriented graph $D$ of minimum semi-degree $d$ of order $n\leq rd$ and of connectivity $k(D)\leq\psi(r)d$, contains a triangle ?

\vspace{0.3cm} By the result of \cite{Lichia}, we have $\psi(r)=1$ for $2<r\leq 2.91082$. By Theorem 1.3, for $2.91082<r<\frac{2}{c}$ we have 
$\psi(r)\geq \max \left\lbrace \dfrac{5-r-4c+c^2}{(1-c)(2-c)}d,\,\dfrac{2-c r}{2-c}d\right\rbrace$.
Thus, since $c\leq 0.3465$, we get  $\psi (3)> 0.679$, $\psi (3.5)> 0.476$, $\psi (4) > 0.371$ $\psi(4.5)>0.266$, $\psi (5)> 0.161$ and $\psi (5.5)> 0.057$. Observe that Conjecture 1.1 is true, if and only if $\psi (3)=1$.



\begin{thebibliography}{}
\bibitem{Bang} J. Bang-Jensen and G. Gutin, Digraphs, Springer,
2002, p. 555.
\bibitem {Li} H. J Broersma and X. Li, Some approaches to a conjecture on short cycles in Digraphs, {\it Discrete Applied Math}, 120 (1-3), (2002), 45-53.
\bibitem{Caccetta} L. Caccetta and R. H\"{a}ggkvist, On minimal
digraphs with given girth, {\it Congressus Numerantium}, 21, 1978, p.
181-187.
\bibitem{Chen} K. Chen, S. Karson, D. Liu and J. Shen, On the
Chudnovski-Seymour-Sullivan conjecture on cycles in triangle-free
digraphs, arXiv:0909.2468v1[math.CO], submitted to Discrete Math.
\bibitem{Chud} M. Chudnovsky, P. Seymour and B. Sullivan, Cycles in
dense digraphs, Combinatorica 28 (2008), 1-18
\bibitem{Hamb} P. Hamburger, P. Haxell and A. Kostochka, On directed triangles in
digraphs, Electronic Journ. of Combinatorics., 14 (2007) N19.
\bibitem{Kral} J. Hladk\'{y}, D. Kr\'{a}l', S. Norin, Counting flags in triangle-free digraphs, {\it Electronic Notes in Discrete Math }, 34 (2009), 621-625.
\bibitem {Lichia} N. Lichiardopol, A new bound for a particular case of the Caccetta-H\"{a}ggkvist conjecture, {\it Discrete Math.}, 310(23) (2010), 3368-3372.
\bibitem{Lichiard} N. Lichiardopol, A new lower bound on the strong connectivity of an oriented graph. Application to diameters with a particular case related to Caccetta-H\"{a}ggkvist conjecture {\it Discrete Math}, {308(22)} (2008), 5274-5279.
\bibitem{Kuhn} L. Kelly, D. K\"{u}hn and D. Osthus, Cycles of given length in oriented graphs, {\it Journal of Comb. Theory, Series B}, 100 (3) (2010), 251-264.



\end{thebibliography}
\end{document}